\documentclass[draft]{amsart}

\usepackage[english]{babel}

\usepackage{amssymb}
\usepackage{amsmath}
\usepackage{amsthm}
\usepackage{enumerate}
\usepackage{cite}
\usepackage{xcolor}
\usepackage[final]{listings}
\usepackage[section]{algorithm}

\title[Hilbert functions and Tensor Analysis]
{Hilbert functions and Tensor Analysis}
\date{}

\newcommand{\C}{\mathbb{C}}
\newcommand{\N}{\mathbb{N}}
\newcommand{\Z}{\mathbb{Z}}
\newcommand{\Pj}{\mathbb{P}}

\newcommand{\vect}[1]{\mathbf{#1}}

\newcommand{\Tang}[2]{\mathrm{T}_{#1} {#2}}

\newcommand{\rank}{\operatorname{rank}}

\newtheorem{defn0}{Definition}[section]
\newtheorem{prop0}[defn0]{Proposition}
\newtheorem{thm0}[defn0]{Theorem}
\newtheorem{lemma0}[defn0]{Lemma}

\newtheorem{exa0}[defn0]{Example}
\newtheorem{rem0}[defn0]{Remark}


\author[L.~Chiantini]{Luca Chiantini}
\address{Dipartimento di Ingegneria dell'Informazione e Scienze Matematiche, Universit\`a di Siena, Italy}
\email{luca.chiantini@unisi.it}

\begin{document}

\begin{abstract} We show how well known tools of algebraic geometry for the study of finite sets
can be fruitfully applied to the study of Waring decompositions of symmetric tensors (forms). We mainly focus on the 
uniqueness of a given decomposition (the identifiability problem), and show
how, in some cases, one can effectively determine the uniqueness even in some range in which
the Kruskal's criterion does not apply. 
\end{abstract}

\maketitle

\section{Introduction }

The paper aims to introduce some basic geometric methods for the study of the
decompositions of tensors. It is mainly devoted to symmetric decompositions of symmetric tensors, which 
can be identified with homogeneous polynomials, i.e. forms.

Decomposing a form $F$ as a sum of powers (Waring decomposition) is a crucial step to understand the
complexity of $F$. The complexity, or (Waring) rank, of $F$ is indeed given by the minimal number of summands 
which are necessary to express $F$ as a sum of powers.

In many effective cases, it turns out that one has one decomposition of $F$ as a sum of powers, and the
problem is   to determine if the given decomposition has minimal length or it is unique (up to trivialities). Just to give a 
couple of examples:

- in the Strassen problem, one has a form which is  a sum $F=F_1+F_2$ where $F_1,F_2$ are forms defined
over two different, disjoint sets of variables. Then one can assume to have a minimal decomposition of both $F_1$
and $F_2$. The problem is to determine if the sum of the two decompositions gives a  decomposition
of $F$ {\it of minimal length}.  See \cite{CarlCatalC15}  and \cite{Shitovb}, for recent accounts on the theory.

- in the application of tensor analysis to signal processing, there are computational methods which can determine
(an approximation of) one decomposition of a tensor $F$. Since one aims to reconstruct the original components
of a mixed signal, the uniqueness of the decomposition is crucial to guarantee that the computed decomposition
is (in a small neighborhood of) the correct one (see e.g. \cite{RaoLiZhang18}).

For the identifiability problem, i.e. in order to determine that a decomposition is unique (up to trivialities), the most popular
criterion  is the Kruskal's criterion (see Theorem \ref{kruthm} below), which requires the calculation of
the Kruskal rank of a set of points (see Definition \ref{krudef}). Kruskal's criterion only works for small values of the rank.
Recently, for symmetric tensors, there is a series of results which show how the Kruskal's criterion can be
modified, to widen slightly the range of application (see \cite{BallBern12a}, \cite{BallC12}, \cite{COttVan17b}, 
\cite{AngeCVan}). These extensions of Kruskal's criterion are mainly based
on methods of algebraic geometry for the study of finite sets in projective spaces.

Since we believe that geometric tools for the study of finite projective sets can contribute to many other
aspects of the theory of symmetric (and maybe also non-symmetric) tensors, and we feel that several
tools are not widely known in the community of researchers in tensor analysis, we provide here an account of
methods which constitute the background for the theory
developed in \cite{COttVan17b} and \cite{AngeCVan}.

As a by-product, we show how similar argument yield a slight extension of the results of \cite{AngeCVan}, for
forms of degree $4$, even
to the case in which the Kruskal rank of a given decomposition is not maximal (see Theorem \ref{quartplus}).

We hope, in this way, to contribute to the propagation of geometric tools which can help a lot
 our insight into the analysis of decompositions of specific tensors. 
 
The structure of the paper is the following.
The first section contains some basic definitions, basic results and remarks which are useful in the theory.
The second section contains a list of results on tensors which are proved by means of the Hilbert function.
The third section is devoted to prove a new result, which extends a recent criterion, proved
by Angelini, Vannieuwenhoven and the author (\cite{AngeCVan}), for the (symmetric) identifiability
of a symmetric tensor in a range where the  Kruskal's criterion does not apply. 
The  result requires a deep analysis of the Hilbert function of  a finite set in a projective space.
In the last section there is a short list of possible developments of the theory and open problems.

\section{Tensors and projective geometry}

Since the study of tensors under a geometric point of view is strictly related with systems of homogeneous
polynomials and their solutions, it is natural, from a mathematical point of view, to treat tensors defined over 
an algebraically closed field, as the complex field $\C$. 

At the risks of losing a strict connection with experience, yet the choice of working over $\C$ will not
sound so odd to specialists of quantum information theory, where the algebraic properties of complex
numbers play a primary role in many quantum manipulations.

Less familiar is the choice of working on \emph{projective} spaces of tensors. The idea behind using the 
projective setting is that the phenomena encoded in a tensor $T$ are as well encoded
in its multiples $aT$, for $a\in\C$ a non-zero constant. In projective spaces, a point $P$ is an equivalence
class containing a vector and its multiples. At the cost of dropping the one-to-one correspondence between
points and coordinates (which are defined up to \emph{scaling}), projective geometry provides a compact
algebraic ambient where some operations, like linear dependence, have a natural interpretation.

Thus, we drop the {\it probabilistic} approach, in which the sum of some entries of the tensors
are forced to be $1$, since they represent the probabilities of some event,  and we will freely multiply tensors by (complex) scalars.
It is an ubiquitous fact that all the results that we obtain can be translated in the probabilistic language,
without any loss of validity. The main, non-trivial aspect of the projective point of view 
is the notion of \emph{product} of projective spaces, which does not produce a linear variety. 
\bigskip

So, we consider a complex vector space $V$ of dimension $n+1$, which we will often identify with $\C^{n+1}$, thanks 
to the choice of a basis. We will think of $V$ as the space of \emph{linear forms} $a_0x_0+a_1x_1+\dots+a_nx_n$,
where $x_1,\dots, x_n$ can be identified with the elements of the chosen  basis or with variables.
Consequently, the space $Sym^d(V)=Sym^d(\C^{n+1})$ will be
identified with the space of homogeneous polynomials (forms) of degree $d$ in the $n+1$ variables $x_0,\dots, x_n$.

Instead of considering directly symmetric tensors as vectors of $Sym^d(V)$, we consider the projective space
$\Pj(Sym^d(V))$ and consider points $T$ in this space. Thus $T$ corresponds to a symmetric tensor or
a form, modulo scaling. Any representative for the equivalence of class of $T$ is a \emph{set
of coordinates} for $T$. As $Sym^d(V)$ has dimension $\binom{n+d}d$, the space $\Pj(Sym^d(V))$
has projective dimension 
$$N(d,n):=\binom{n+d}d-1.$$
\smallskip

The next step is the definition of a (non-linear) map from $\Pj(V)=\Pj^n$ to the space $\Pj(Sym^d(V))=\Pj^{N(d,n)}$:
the Veronese map. 

To do that, choose an order for the monomials of degree $d$ in $n+1$
variables $M_0,\dots, M_N$, $N=N(d,n)$. One of the most popular order is the lexicographic one,
and we will opt for it for the rest of the paper. 

Then, use the  coordinates to define a map $\nu_{d,n}$ as follows.
Let a point $P\in\Pj(V)$ have coordinates $a_0x_0+\dots+a_nx_n$.
We will write:
$$P=[a_0x_0+\dots+a_nx_n].$$ 
We define $\nu_{d,n}$ by sending $P$ to the equivalence class
$$\nu_{d,n}(P)=[(a_0x_0+\dots+a_nx_n)^d].$$
The class $\nu_{d,n}(P)$  does not depend on the choice of a representative for the class $P$, so we get 
a well defined projective map. We will refer to this map as the \emph{Veronese map of degree $d$ in $n+1$ variables}.
We will often write $\nu_d$ for the Veronese map, when there is no confusion on the number of variables.

We notice that the Veronese maps are embeddings. 

\begin{prop0} Every Veronese map  $\nu_{d,n}$ is injective.
\end{prop0}
\begin{proof} Assume that two points $P,Q\in\Pj^n$ have the same image in $\nu_{d,n}$.
Choose coordinates in $\Pj(V)$ and let $a_0x_0+\dots+a_nx_n$ be coordinates for $P$ and
$b_0x_0+\dots+b_nx_n$ be coordinates for $Q$. Then $(b_0x_0+\dots+b_nx_n)^d$ is 
equal to $\alpha(a_0x_0+\dots+a_nx_n)^d$, for some $\alpha\in\C\setminus\{0\}$. Since $\C$ is
algebraically closed, then, after scaling $a_0x_0+\dots+a_nx_n$ by a $d$-root of $\alpha$, 
we may assume $(b_0x_0+\dots+b_nx_n)^d=(a_0x_0+\dots+a_nx_n)^d$. Thus
$b_i=\epsilon_ia_i$, for some choice of the $d$-roots of unit $\epsilon_i$, $i=0,\dots,n$. We want to
prove that the $\epsilon_i$'s are all equal, so that $P=Q$. Indeed, since $\epsilon_0^{(d-j)}\epsilon_i^j=1$ for all $i,j$,
multiplying by $\epsilon_0^j$ it follows $\epsilon_0^j=\epsilon_i^j$ for any $j$, hence $\epsilon_0=\epsilon_i$
for all $i$.
\end{proof}

Notice that the previous construction is not the unique way to define a Veronese map. Often $v_{d,n}(P)$ is defined by
computing $b_i=M_i(a_0,\dots,a_n)$ for $i=0,\dots,N$ and sending $P$ to the equivalence class
$[b_0M_0+\dots +b_NM_N].$
We made our choice in order to make it obvious that the image of the Veronese map is the set of 
forms which are a power of a linear forms. Since the two choices differ only by the multiplication by a
non-singular diagonal matrix, the geometric properties will not be affected after taking any of the choices.
\medskip

Next, we  need to fix some notation for finite subsets of a projective space.

Let $A\subset \Pj^n$ be a non-empty finite set. We denote by $\ell(A)$ the cardinality of $A$. We will say that $A$ is
\emph{linearly independent} when choosing a set of coordinates for each point of $A$ we
get a set of linearly independent vectors. This definition does not depend on the choice
of the coordinates for each point.

We will denote with $\langle A\rangle$ the \emph{linear span} of $A$. 

\begin{rem0}\label{span} {\rm
The projective dimension of $\langle A\rangle$
is at most $\ell(A)-1$. The dimension of $\langle A\rangle$ is equal to $\ell(A)-1$ precisely when $A$ is linearly independent. 

Notice that, by elementary linear algebra, for any finite set $A\subset \Pj^n=\Pj(V)$ the dimension of the linear span $\langle A\rangle$
is equal to $n$ minus the dimension of the space of linear forms that vanish
at the points of $A$.
}\end{rem0}

\begin{defn0}\label{krudef}
Let $ A \subset \Pj^n$ be a finite set. The \emph{Kruskal rank}  is the maximum integer $k_A$
such that any subset $B\subset A$ of cardinality $\ell(B)\leq k_A$ is linearly independent. 
\end{defn0}

Notice that $k_A$ is at most equal to $\ell(A)$, and $k_A=\ell(A)$ if and only if $A$ is linearly independent. 
Unless $A$ is a singleton, then $k_A$ is always bigger than $1$. Moreover $k_A=2$ exactly when $A$ is aligned.

Obviously the Kruskal rank of a set of points $A\subset\Pj^n$ cannot exceed neither $n+1$, nor the cardinality of $A$.
We have indeed:
$$ k_A\leq \dim\langle A \rangle + 1 \leq \ell(A).$$
Next definition concerns the case where the Kruskal rank is maximal.

\begin{defn0}\label{LGP}
A finite set $ A \subset \Pj^n $ is in \emph{linear general position} (LGP) if the Kruskal rank of $A$ is maximal, i.e.
the Kruskal rank is equal to $\min\{\ell(A), n+1\}$.
This is equivalent to say that for any $a\leq n+1$, any subset of $A$ of cardinality $a$ is linearly independent. 
\end{defn0}

Next, we come to the definition of \emph{decomposition} of a (symmetric) tensor.

\begin{defn0}
Let $A \subset \Pj^n=\Pj(V)$ be a finite set.  We say that $A$ is a \emph{decomposition} of the tensor  $T\in \Pj(Sym^d(V))$, or equivalently
that $A$ \emph{computes} $T$, if $T$ belongs to the span $\langle \nu_{d}(A) \rangle$.
\end{defn0}

\begin{defn0}
Let $ A \subset \Pj^n $ be a decomposition of  $T$. $ A $ is \emph{minimal} if we cannot find a proper subset $ A' $ of $ A $ 
such that $T \in \langle\nu_{d}(A')\rangle $.
\end{defn0}

\begin{rem0}\label{rem:indep}{\rm
If $ A \subset \Pj^n $ is a decomposition of $T$ and satisfies the minimality property, then in particular 
the points of $ \nu_{d}(A) $ are  linearly independent, i.e., $$ \dim(\langle\nu_{d}(A)\rangle) = \ell(A) -1.  $$}
\end{rem0}

\subsection{The Hilbert function of finite sets in projective spaces}\label{sec:hilb}

We collect in this section a series of definitions and propositions which are well known to people working 
in algebraic geometry, but maybe not so familiar to other people working in tensor analysis. The main definition is the 
\emph{Hilbert function} of a finite set in a projective space, which is a basic tool for our results on the decompositions
of symmetric tensors.

\begin{defn0}
Let $Y\subset \C^{n+1} $ be an ordered, finite set of cardinality $\ell $ of vectors. Fix an integer $ d \in \N $. 

The \emph{evaluation map of degree $d$ on $Y$} is the linear map
$$ ev_{Y}(d): Sym^d(\C^{n+1}) \to \C^\ell $$ 
which sends $ F \in Sym^d(\C^{n+1}) $ to the evaluation of $ F$ at the vectors of $Y$. 
\end{defn0}

We will use the evaluation map to define the Hilbert function of a finite set $ Z \subset \Pj^n $.

\begin{rem0}\label{otherdef}{\rm
Let $A \subset \Pj^n $ be a finite set, with a definite order. Choose a set of homogeneous coordinates for the points of $A$.
We get an ordered set of vectors $Y\subset \C^{n+1} $, for which the evaluation map  $ev_{Y}(d)$  is defined for every $d$.

If we change the choice of the homogeneous coordinates for the points of the fixed set $A$, we get another ordered
set $ Y' \subset \C^{n+1} $ and the evaluation map $ev_{Y'}(j)$ differs from $ev_{Y}(j)$ for the
multiplication by a non-singular diagonal matrix. 
Thus the rank of $ev_{Y}(j)$ and $ev_{Y'}(j)$ are the same for all $j$.

It is also clear that the rank of $ev_{Y}(j)$ does not depend on how we ordered the points of $A$.

Let $f:\C^{n+1}\to\C^{n+1}$ be an automorphism and consider the associated change of coordinates $\Pj^n\to\Pj^n$,
that we call again $f$, by abuse. Then the evaluation on $Y$ and $f(Y)$ differ by the multiplication by
a non-singular matrix. Thus for any $d$ the maps $ev_Y(d)$ and $ev_{f(Y)}(d)$ have the same rank.}
\end{rem0}

\begin{defn0}
Let $Z\subset \Pj^n $ be a finite set. Choose an order and an ordered  set of homogeneous coordinates $Y$  for the points of $A$.
Define the \emph{Hilbert function} of $Z$ as the map
$$ h_Z : \Z \to \N \qquad h_Z(d) = \rank(ev_{Y}(d)) .$$ 
By the previous remark, the Hilbert function does not depend on the choice of the coordinates, as well as it 
does not vary after a change of coordinates in $\Pj^n$.
\end{defn0}

People who are expert of algebraic geometry may wonder why we did not define the Hilbert function
as the rank of the restriction maps $H^0(\mathcal O(d))\to H^0(\mathcal O_Z(d))$, where $\mathcal O,\mathcal O_Z$
indicate respectively the structure sheaves of $\Pj^n$ and $A$. This would simplify the notation,
since the restriction is well defined, regardless of a choice of coordinates for the points of $A$.
On the other hand, our definition is immediately accessible also to readers who are not expert about
cohomology, structure sheaves and so on. We preferred to make our basic definition more familiar and easily 
computable for a wider audience. We based our definition on the choice of coordinates because
only after a choice of coordinates for the points of $A$ one has a natural identification
of $ H^0(\mathcal O_Z(d))$ with $\C^{\ell}$.

There is a different notation for the Hilbert function, which is widely used in algebraic geometry.
Since it  clarifies some aspects, we introduce it.

\begin{rem0}\label{cond}{\rm
Recall that the homogeneous ideal $I_Z$ of the set $Z$ in the polynomial ring $\C[t_0,\dots,t_n]$ 
is the ideal generated by all the homogeneous polynomials (forms) which vanish at all the points of $Z$.
Thus, $I_Z$ is a graded ideal. Its degree $d$ summand $I_Z(d)$ is exactly the kernel of the evaluation map $ev_Z(d)$.

Notice that, indeed, the kernel does not depend on the choice of homogeneous coordinates for the points of $Z$,
because the vanishing of a form at a projective point $P$ is independent from the choice of a specific
set of homogeneous coordinates for $P$. 

Thus, recalling that the vector space of forms of degree $d$  we have
$$h_Z(d)= \dim(Sym^d(\C^{n+1}))- \dim(I_Z(d))=\binom{n+d}n-\dim(I_Z(d)).$$}
\end{rem0}

Consequently, we introduce the following notation:

\begin{defn0}\label{conds}
Let $Z$ be a finite subset of the projective space $\Pj^n$and let $h_Z$ be its Hilbert function.
 For any $d \geq 0$, the value $ h_Z(d) $ is also called the \emph{number of conditions 
that $Z$ imposes to  forms of degree $d$}.

We say that $Z$ \emph{imposes independent conditions to forms of degree $d$}, or also that 
\emph{the points of $Z$ are separated by forms of degree $d$}, if $ h_Z(j) = \ell(Z) $. 
This happens exactly when, for (any choice of) a set $Y$  of homogeneous coordinates  for the points of $Z$,
the evaluation map $ev_Y(d)$ surjects.
\end{defn0}

\begin{rem0}{\rm Let us explain in more details the last definition. Set $\ell=\ell(Z)$, and fix an order for the points
of $Z$. 

Take a vector $e_j=(0,\dots,0,1,0,\dots,0)$ ($1$ is in the $j$-th position) of the natural basis of $\C^\ell$, 
which corresponds to the $j$-th point $P_j$ of $Z$ in the given order. 
We say that \emph{$P_j$ is separated in $Z$ by forms of degree $d$} if $e_j$ belongs to the image of the
evaluation map $ev_Y(d)$. Indeed, in this case, $e_j$ is the evaluation of
a form $F$ of degree $d$. Thus there exists a form $F$ which vanishes at all the points of $Z$, but $P_j$.
Notice that this is independent on the choice of the homogeneous coordinates $Y$.

If $ h_Z(j) = \ell(Z) $, i.e. if the evaluation map $ev_Y(d)$ surjects, then any point of $Z$ is separated.}
\end{rem0}

The link between the Hilbert function of finite sets and the decompositions of symmetric tensors is mainly based on the following formula,
which gives a different, geometric interpretation of the values $h_Z(d) $.

\begin{prop0}\label{dimspan}
Let $\nu_{d,n}: \Pj^n\to \Pj^N$, $N=N(d,n)$, be the $d$-th Veronese embedding of $\Pj^n$.
For any finite set $Z\subset \Pj^n$, and for any $d \geq 0 $, the value $ h_Z(d) $ determines the dimension 
of the span of $\nu_d(Z)$. I.e.:
$$ h_Z(d) = \dim (\langle \nu_{d,n}(Z) \rangle) +1. $$
\end{prop0}
\begin{proof} We know that the value $h_Z(d)$ is equal to the dimension of $Sym^d(\C^{n+1})$ minus the
dimension of the space $I_Z(d)$, where $I_Z$ is the homogeneous ideal of $Z$ in $\C[t_0,\dots,t_n]$.
If we identify the  coordinates in $\Pj^{N(d,n)}=\Pj(Sym^d(\C^{n+1}))$ with the monic monomials $M_j$'s
of degree $d$ in $\C[t_0,\dots,t_n]$ (say with the lexicographic order), then any element of $I_Z(d)$
corresponds to a linear form in $\Pj(Sym^d(\C^{n+1}))$. The claim follows by Remark \ref{span}.
\end{proof}

\begin{defn0}
We define the \emph{first difference of the Hilbert function} $Dh_Z$ of $Z$ as:
$$ Dh_Z(j) = h_Z(j)-h_Z(j-1),\quad j \in \Z .$$
The set of non-zero values of $Dh_Z$ is called the \emph{h-vector} of $Z$.
\end{defn0}

The following properties of $h_A$ and $Dh_A$ are elementary and well-known in algebraic geometry. 
We recall them because they  will be useful throughout the paper.

\begin{lemma0}\label{rem:triv} Set $\ell=\ell(Z)$. Then we have:
\begin{enumerate}[(i)]
\item $ h_Z(d) \leq \ell$ for all $d$;
\item $ Dh_Z(d) = 0$ for $d < 0$;
\item $ h_Z(0) = Dh_Z(0) = 1$;
\item $ Dh_Z(d) \geq 0 $ for all $d$;
\item $ h_Z(d) = \ell(Z) $ for all $d \geq \ell(Z)-1$; 
\item $ h_Z(i) = \sum_{0\leq d \leq i} Dh_Z(d) $; 
\item $ Dh_Z(d) = 0$ for $ d \gg 0 $ and $ \sum_{d} Dh_Z(d) = \ell(Z) $; 
\item if $ h_Z(d) = \ell(Z) $, then $ Dh_Z(d+1) = 0 $.
\end{enumerate}
\end{lemma0}
\begin{proof} (i) is a consequence of the definition. (ii)  follows immediately since the space $Sym^d(\C^{n+1})$ is $(0)$ for $d$ negative. 
(iii) follows since $Sym^0(\C^{n+1})=\C$ and the evaluation of a constant  form $c$ is equal to $c(1,\dots,1)\in\C^\ell$.

To see (iv), fix an ordered set of coordinates $Y$ for the points of $Z$ and fix a linear form $\Lambda$ which does not vanish at
any vector of the finite set $Y$. Then for any form $F$ of degree $d$, the evaluation of $\Lambda F$ at $Y$ is equal to the
evaluation of $F$ at $Y$ multiplied by a fixed non-singular diagonal matrix, whose entries are the evaluations of $\Lambda$
at the vectors of $Y$. Thus the image of $ev_Y(d+1)$ contains a subspace isomorphic to the image of $ev_Y(d)$.
It follows that $h_Z(d+1)\geq h_Z(d)$, hence $Dh_Z(d)\geq 0$.

To see (v), choose for each point $P_j\in Z$ a linear form $L_j$ which vanishes at $P_j$ and does not vanish
at any other point $P_k\in Z$. Then for any $j$ call $F_j$ the product of the linear forms $L_k$, $k\neq j$.
$F_j$ is a form of degree $\ell-1$, which vanishes at all the points of $Z$, except $P_j$. Thus,
the evaluation of $F_j$ at an ordered set of coordinates $Y$ for the points of $Z$ is a vector
$(c_1,\dots,c_\ell)$ with $c_k=0$ for $k\neq j$ and $c_j\neq 0$. It follows that $ev_Y(\ell-1)$ is surjective.
Then, by (iv), $ev_Y(d)$ surjects for all $d\geq \ell-1$.

(vi) is a triviality. (vii) and (viii) are obvious consequences of (v) and (vi).
\end{proof}

\begin{prop0} \label{inclu} With the previous notation, if $ Z' \subset Z$, then, for every $d \in \Z $, 
we have $h_{Z'}(d) \leq h_Z(d) $ and $Dh_{Z'}(d) \leq Dh_Z(d).$
\end{prop0}
\begin{proof}
Fix, as usual, an ordered set of coordinates $Y',Y$ for the points of $Z',Z$ respectively. Then we have an obvious forgetful map
$f:C^\ell\to\C^{\ell'}$, where $\ell'=\ell(Z')$, such that $ev_{Y'}(d)=f\circ ev_Y(d)$ for all $d$. This implies 
that $h_{Z'}(d) \leq h_Z(d) $.

The second inequality is less trivial, and we will need some algebra. 
Write $R$ for the polynomial ring $\C[t_0,\dots,t_n]$ and call $I_Z$ the ideal generated by forms which vanish at
the points of $Z$. The inclusion $I_Y\subset R$ determines, for every $d\in\Z$ an exact sequence of vector spaces:
$$ 0\to I_Y(d)\to R(d)\to (R/I)(d)\to 0,$$
where $R(d), R/I_Z(d)$ are the graded pieces of the rings $R,R/I$ respectively, in degree $d$.
It follows by  Remark \ref{otherdef} that for any $d$:
$$h_Z(d)=\dim(R/I_Z(d)).$$
The natural inclusion $I_Z\subset I_{Z'}$ induces a surjection $R/I_Z(d)\to R/I_{Z'}(d)$ for all $d$.
Let $\Lambda$ be a linear form  in $\C[t_0,\dots,t_n]$, which does not vanish at any point of $Z$.
The multiplication by $\Lambda$ induces an inclusion $R/I_Z(d)\to R/I_Z(d+1)$. Indeed if $F\in R(d)$ is a form which does not
vanish at some point $P\in Z$, then $LF$ cannot vanish at $P$, i.e. the class of $LF$ is non-zero in $R/I_Z(d+1)$.
Call $J_Z$ the ideal generated by $I_Z$ and $\Lambda$. We have an exact sequence:
$$0\to R/I_Z(d)\to R/I_Z(d+1)\to R/J_Z(d+1)\to 0$$
which proves that
$$Dh_Z(d)= \dim(R/J_Z(d+1)).$$
Similarly $\Lambda$ induces an embedding $R/I_{Z'}(d)\to R/I_{Z'}(d+1)$ and $Dh_{Z'}(d)= \dim(R/J_{Z'}(d)).$ Now look at
the commutative diagram:
$$\begin{matrix}
0 & \to & R/I_Z(d)    &\stackrel{L}\longrightarrow & R/I_Z(d+1) & \to & R/J_Z(d+1)  &\to &0 \\
  &  &  \downarrow  &                      & \downarrow  &      & \downarrow & & \\
0 & \to & R/I_{Z'}(d)  &\stackrel{L}\longrightarrow & R/I_{Z'}(d+1) & \to & R/J_{Z'}(d+1)  &\to &0 
\end{matrix}
$$
Since the central vertical map $R/I_Z(d+1)\to R/I_{Z'}(d+1)$ surjects, by the snake's lemma also the map
$R/J_Z(d+1)\to R/J_{Z'}(d+1)$ surjects. Then $Dh_Z(d)=\dim( R/J_Z(d+1))\geq \dim(R/J_{Z'}(d+1))=Dh_{Z'}(d)$. 
This proves the second claim.
\end{proof}

Perhaps, the most important algebraic result on Hilbert functions of finite sets is the {\it maximal growth principle}
found by Macaulay.
Roughly speaking, the maximal growth principle gives an upper bound for the value of $h_A(i+1)$ in terms
of $h_Z(i)$ and the dimension of the ambient space. We list below the most relevant consequences
for the application to the study of tensors and forms. 

\begin{prop0}\label{nonincr}
Assume that  for some $j>0$ we have $Dh_Z(j) \leq j$. Then:
$$ Dh_Z(j) \geq Dh_Z(j+1). $$
In particular, if  for some $j>0$, $ Dh_Z(j)=0 $, then $Dh_Z(i)=0$ for all $i\geq j$.
\end{prop0}
\begin{proof} See section $3$ of  \cite{BigaGerMig94}.
\end{proof}

\begin{exa0} \label{Exh(1)} {\rm Let us see what happens for $h_Z(1)$. Since for $i=1$ the domain of the evaluation map
is $Sym^1(\C^{n+1})=\C^{n+1}$, then clearly $h_Z(1)\leq n+1$. So $h_Z(1)=0$ can hold only if
$\ell(Z)\leq n+1$. Moreover the kernel of the evaluation map  $ev_Z(1)$ is isomorphic to the space
of linear forms in $\Pj^n$ which vanish at $Z$. Thus:
$$ h_Z(1)= 1+\dim(\langle Z \rangle).$$
In particular, $ h_Z(1)=0$ if and only if $Z$ is linearly independent.}
\end{exa0}

\begin{rem0}\label{indepcond} {\rm
Assume that for some $j$ we have  $ Dh_Z(j) = 0 $, so that $ h_Z(j-1) = h_Z(j) $. By Proposition \ref{nonincr},
 for any $ i \geq j $ also $ Dh_Z(i) = 0 $, i.e.,
$ h_Z(j-1)=h_Z(i)$ for any $ i \geq j $. Therefore, by part (v)  of Lemma \ref{rem:triv}, 
 $h_Z(j-1)$ is equal to the cardinality of $Z$, i.e., the evaluation map in degree $j-1$ surjects and $Z$ imposes
 independent conditions to hypersurfaces of degree $j-1$.}
\end{rem0}

\begin{rem0}\label{resto1}{\rm
Assume $h_Z(i)=\ell(Z)-1$. Then $h_Z(i+1)>h_Z(i)$, by Remark \ref{indepcond}.
Thus, if $h_Z(i)=\ell(Z)-1$, then necessarily $h_Z(i+1)=\ell(Z)$.}
\end{rem0}

Hilbert functions of finite sets share many other properties.
One can find an accurate account of the theory in the book of Iarrobino and Kanev \cite{IK} and
in the book of Migliore  \cite{Migliore}.

We will focus on the {\it Cayley-Bacharach} property, which is defined as follows:

\begin{defn0}\label{def:CB}
A finite set $Z\subset \Pj^n$ satisfies the \emph{Cayley-Bacharach property in degree $i$}, 
abbreviated as $\mathit{CB}(i)$, if, for any $P \in Z$, 
 every form of degree $i$ vanishing at $ Z\setminus\{ P\}$
also vanishes at $P$.
\end{defn0}

\begin{rem0}{\rm One should compare $\mathit{CB}$ with the property of separating points. In a sort of sense, the $\mathit{CB}$ property is
the contrary of the separation property. 

- $Z$ is separated in degree $i$ if for all $P \in Z$, there exists a form of degree $i$ vanishing at $ Z\setminus\{ P\}$
and not vanishing at $P$.

- $Z$ does not satisfy $\mathit{CB}$ if there exists $P\in Z$ and there exists a form of degree $i$ vanishing at $ Z\setminus\{ P\}$
and not vanishing at $P$.

In particular, if $Z$ satisfies $\mathit CB(i)$,  then  hypersurfaces of degree $i$ cannot separate the points of $Z$,
i.e. $h_Z(i)<\ell(Z)$.
} \end{rem0}

\begin{exa0} {\rm
The set $ Z $ consisting of four points in $ \Pj^2 $, three of them aligned,
does not satisfy $\mathit{CB}(1)$, and $h_{Z}(1)<4$.

Let $ Z$ be a set of $6$ points in $\Pj^{2} $. 

If the $6$ points are general, then $ Dh_{Z} = (1,2,3) $,  and $ Z $ satisfies $\mathit{CB}(1)$.
Since $h_Z(2)=6$, $Z$ does not satisfy $CB(2)$.

If $Z$ lies on an irreducible conic, then $ Dh_{Z} = (1,2,2,1) $,  
 and $ Z $ satisfies $\mathit{CB}(2)$, and, hence, $\mathit{CB}(1)$.
 
If $Z$ has $ 5 $ points on a line plus one point off the line, then 
 $ Dh_{Z} = (1,2,1,1,1) $,  and $ Z $ does not satisfy $\mathit{CB}(1)$.}
\end{exa0}

\begin{rem0}\label{rem:CBprop}{\rm
If $Z$ satisfies $\mathit{CB}(i)$, then it satisfies $\mathit{CB}(i-1)$ too. Otherwise, one could find 
$ P \in Z $ and a  hypersurface $ F \subset \Pj^{n} $ of degree $ (i-1) $ such that $ Z \setminus \{P\} 
\subset F $ and $ P \notin F $. Therefore, if $ H_{P} \subset \Pj^{n} $ is a hyperplane 
not containing $ P $, then $ F \cup H_{P} \in H^{0}(J_{Z\setminus\{ P\}}(i)) \setminus H^{0}(J_{Z}(i)) $, 
which contradicts the hypothesis.}
\end{rem0}

\begin{rem0} {\rm
Assume that $Z$ satisfies $\mathit{CB}(i)$. Call $I_Z$ the homogeneous ideal of $Z$. For any $ P \in Z $ call 
$I_{Z\setminus\{ P\}}$ the homogeneous ideal of $Z\setminus\{ P\}$. Then for all $j\leq i$ we have
$ I_Z$ and $I_{Z\setminus\{ P\}}$ are equal in degree $j$. It follows that:
\begin{equation}\label{eq:h0}
h_Z(j)=h_{Z\setminus\{P\}}(j) \quad \mbox{ and }\quad Dh_Z(j)=Dh_{Z\setminus\{P\}}(j) \quad \forall j\leq i.
\end{equation}
}
\end{rem0}

The following proposition, which gives a strong bound on the Hilbert function of sets with a Cayley-Bacharach property,
 is a refinement of a result due to Geramita, Kreuzer, and Robbiano (see Corollary 3.7 part (b) and (c) of \cite{GerKreuzerRobbiano93}).

\begin{thm0}\label{GKRext}
If a finite set $ Z \subset \Pj^{n} $ satisfies $\mathit{CB}(i)$, then for any $ j $ such that $ 0 \leq j \leq i+1 $ we have
$$ Dh_{Z}(0)+Dh_{Z}(1)+\cdots + Dh_{Z}(j) \leq Dh_{Z}(i+1-j)+\cdots +Dh_{Z}(i+1).$$
\end{thm0}

\begin{proof} 
See Theorem 4.9 of \cite{AngeCVan}.
\end{proof}

Finally, let us point out the relation between the Hilbert functions of a finite set $Z$ and of its image
in a Veronese map $\nu_d(Z)$.

\begin{rem0} \label{remh(1)}{\rm  Let $Z\subset\Pj^n$ be a finite set and let $\nu_d(Z)\subset\Pj^N$ be its image
in the $d$-th Veronese map. Then
$$ h_Z(d)=h_{\nu_d(Z)}(1).$$
Namely the inverse image in $\nu_d$ of a linear form $\Lambda$ in $\Pj^N$ corresponds to a form of degree
$d$ in $\Pj^n$, and the consequent map $\C^{N+1}\to Sym^d(\C^{n+1})$ surjects. Moreover it is easy to see that,
for any choice of coordinates $Y$ for the points of $Z$ in $\Pj^n$ and the consequent choice $\nu_d(Y)$ of coordinates
for the points of $\nu_d(Z)$, one has $ev_{Y'}(L)= ev_Y(\nu_d^{-1}(L))$, so that the claim follows.

In particular, since $\nu_d$ is a bijection, then $h_Z(d)=\ell(Z)$ if and only if $h_{\nu_d(Z)}(1)=\ell(\nu_d(Z))$,
i.e. if and only if $\nu_d(Z)$ is linearly independent (see Example \ref{Exh(1)}).}
\end{rem0}

The following result will be useful in the proof of Theorem \ref{quartplus}

\begin{prop0}\label{v2} Let $Z$ be a finite set in $\Pj^n$. Call $k$ the Kruskal rank of $Z$.
If $\ell(Z)\leq 2k-1$, then $Z$ is separated by forms of degree $2$. Hence  $v_2(Z)$ is 
linearly independent.
\end{prop0}
\begin{proof} We know that $k\leq n+1$. For any point $P\in Z$, consider a partition of the 
residue $Z\setminus \{P\}$ in two disjoint sets $Z_1$, $Z_2$, each of cardinality at most $k-1$.
Since $k-1\leq n$, then the span $L_i$ of $Z_i$ has dimension  strictly smaller than $n$.
Moreover, $L_i$ does not contain $P$, for otherwise $Z$ has $k$ linearly dependent points,
which contradicts the assumption on the Kruskal rank of $Z$.
Thus, there are hyperplanes $H_1,H_2$ containing  $Z_1$ and $Z_2$ respectively and both 
missing $P$. The union $Q=H_1\cup H_2$ is a quadric which
misses $P$ and contains the remaining points of $Z$.
\end{proof}

\section{Results on tensors from classical projective geometry}

The section is devoted to list a series of results on tensors whose
proof is based on the study of the Hilbert function of finite sets. In many cases we
omit the proof, or give only a short draft it.
\smallskip

\begin{rem0} {\rm Fix integers $d,n>1$ and consider symmetric tensors in the space $\Pj(Sym^d(\C^{n+1}))$. 
In \cite{AlexHir95} Alexander and Hirschowitz determined the unique value $r_{d,n}$ such that the set of tensors
of rank $r_{d,n}$ is dense in $\Pj(Sym^d(\C^{n+1}))$. It turns out that $r_{d,n}$ coincides with the expected
value, except for a short list of exceptions. 

We will call $r_{d,n}$ the \emph{generic rank}.}
\end{rem0}

\begin{defn0} We say that a tensor $T\in \Pj(Sym^d(\C^{n+1}))$ of rank $r$ is \emph{identifiable} if $T$ has only one
minimal decomposition $A$ with $\ell(A)=r$, up to scaling and permutations of the summands. 
\end{defn0}

Identifiability is a relevant property for tensors for many applications, as explained in the Introduction.

If we fix a \emph{subgeneric} value of the rank $r< r_{d,n}$, then the set of tensors of rank $\leq r$ in  $\Pj(Sym^d(\C^{n+1}))$ 
is irreducible and its general element has rank $r$,
so we can talk about a \emph{general} tensor of rank $r$. For general tensors of rank $r< r_{d,n}$, thanks to the fundamental
preparatory works \cite{AlexHir95}, \cite{CCi06}, and  \cite{Ball05a}, the situation with respect to 
the identifiability property has been completely described in  \cite{COttVan17a}.

\begin{thm0} Let $d,r\ge 2$. The general tensor in $\Pj(Sym^d(\C^{n+1})$
of subgeneric rank $r < r_{d,n}$ is identifiable, unless it is one of the following cases:
\begin{enumerate}
\item $d=2$;
\item $d=6$, $n=2$, and $r=9$;
\item $d=4$, $n=3$, and $r=8$; 
\item $d=3$, $n=5$, and $r=9$.
\end{enumerate}
In the first case there are infinitely many decompositions. In the three last exceptional cases, 
there are exactly two decompositions.
\end{thm0}
\begin{proof} 
See Theorem 1.1 of \cite{COttVan17a}.
\end{proof}

\begin{rem0}{\rm
On the contrary, when $r=r_{d,n}$, there are very few cases in which a general tensor of rank $r$ is identifiable.
The classification has been proved by Galuppi and Mella, see \cite{GaluppiMella}.

When $r >  r_{d,n}$, the situation is less known. It is not even obvious what is the meaning
of \emph{generic tensors}, since the set of tensors of given rank can have many components.

In any case, one expects that a sufficiently general tensor is \emph{not} identifiable, though for
$r>r_{n,d}$ very few things are known.

For the case $r=r_{n,d}$,  the situation is completely described in \cite{AlexHir95}, \cite{Mella06} and
mainly in \cite{GaluppiMella}: there are  many decompositions, unless $d,n$ are included
in a short list of cases.
}
\end{rem0}

Let us turn to the problem of the identifiability of one \emph{specific} given tensor $T\in \Pj(Sym^d(\C^{n+1}))$,
of which we know a minimal decomposition $A\subset\Pj^n=\Pj(\C^{n+1}))$ with $\ell(A)=r$.

Recall that \emph{minimal} means that the set $v_d(A)$ is linearly independent. We \emph{do not} assume that $\ell(A)$ 
is actually the rank of $T$, i.e. we do not know if $T$ has some other decomposition with smaller cardinality.

Let us start recalling the following, classical result of Sylvester, which disposes of the case $n=1$,
the case of \emph{binary forms}:

\begin{thm0}\label{Sy} Assume $n=1$, i.e. consider the space of tensor $\Pj(Sym^d(\C^{2}))$. Then
$r_{2,d}=(d+1)/2$ if $d$ is odd, $r_{2,d}=(d+2)/2$ if $d$ is even. Moreover every tensor of rank $r<r_{2,d}$
is identifiable.
\end{thm0}
\begin{proof}
See \cite{Sylvester86}.
\end{proof}

Indeed, to be precise, when $n=1$ and $d$ is odd, also tensors of rank $r_{2,d}$ are identifiable. See Theorem \ref{linbound} below.

So, we restrict ourselves to the case $n>1$.
\smallskip

The reason why an analysis of the Hilbert functions is relevant for the identifiability property is expressed in the following
lemma, which can be found in \cite{BallBern12a}:

\begin{lemma0}\label{Hilbident}
Consider two different minimal decompositions  $A,B$ of a tensor $T\in\Pj(Sym^d(\C^{n+1}))$. In other words,
we have:
$$T\in\langle \nu_d(A)\rangle \cap \langle \nu_d(B)\rangle.$$
Then if $Z=A\cup B$, we get $h_Z(d)<\ell(Z)$, so that $Dh_Z(d+1)>0$.
\end{lemma0}
\begin{proof} Set $Z=A\cup B$. First assume that $A,B$ are disjoint.  The existence of $T$ implies that $\nu_d(Z)$
is not linearly independent. By Example \ref{Exh(1)}, this implies that linear forms in the space $\Pj^N$ spanned by $\nu_d(\Pj^n)$ 
do not separate the points of $\nu_d(Z)$. By Remark \ref{remh(1)}, this implies that forms of degree $d$ in $\Pj^n$ do not
separate the points of $Z$. The claims follow by part viii) of  Lemma \ref{rem:triv} and Proposition \ref{nonincr}. 

If $A\cap B\neq \emptyset$, define $B'=A\setminus B$, so that $Z$ is the disjoint union of $A$ and $B'$.
By elementary linear algebra, $\langle\nu_d(A)\rangle \cap \langle\nu_d(B)\rangle$ is also spanned by
$\nu_d(A\cap B)$ and $\langle\nu_d(A)\rangle \cap \langle\nu_d(B')\rangle$. By the minimality of $A$,
$T$ cannot belong to the span of $\nu_d(A\cap B)$. Thus $\langle\nu_d(A)\rangle \cap \langle\nu_d(B')\rangle$
is non empty, and the claim follows again, as above, by part viii) of  Lemma \ref{rem:triv} and Proposition \ref{nonincr}. 
\end{proof}

We can be more precise about the dimension of the intersection of the span of $\nu_d(A)$ and $\nu_d(B)$.

\begin{lemma0}\label{inters} Let $A,B\subset\Pj^n$ be \emph{disjoint} finite sets. Set $Z=A\cup B$.
Then:
$$
\dim(\langle \nu_d(A)\rangle\cap\langle\nu_d(B)\rangle) = \ell(Z)-h_Z(d)-1.
$$
If $A\cap B\neq \emptyset$, then:
$$\dim(\langle \nu_d(A)\rangle\cap\langle\nu_d(B)\rangle)\leq \dim(\nu_d(A\cap B))+\ell(Z)-h_Z(d).$$
\end{lemma0}
\begin{proof} The first formula in an exercise for the application of the Grassmann intersection formula. The second formula
follows since, setting $B_0=B\setminus A$ so that $A,B_0$ are disjoint and $Z=A\cup B_0$, by elementary linear algebra 
$\langle \nu_d(A)\rangle\cap\langle\nu_d(B)\rangle$ is spanned 
by $\nu_d(A\cap B))$ and $\langle \nu_d(A)\rangle\cap\langle\nu_d(B_0)\rangle$.
\end{proof}

An extension of Sylvester's theorem, which works for \emph{all} symmetric tensors in  $\Pj(Sym^d(\C^{n+1}))$,
is possible for $n>1$ only for small values of the rank. The following statement is proved in Theorem 1.5.1
of \cite{BucGinenskyLand13}. We give here an alternative proof, in terms of the Hilbert function of decompositions.

\begin{thm0}\label{linbound} Assume that a tensor $T\in\Pj(Sym^d(\C^{n+1}))$ has a decomposition $A$ with $\ell(A)\leq (d+1)/2$.
Then $T$ has rank $\ell(A)$ and it is identifiable.
\end{thm0}
\begin{proof} Assume on the contrary that $T$ has a second decomposition $B$ with $\ell(B)\leq \ell(A)$, and
take the union $Z=A\cup B$. Then $\ell(Z)\leq 2\ell(A)\leq d+1$. By Lemma \ref{Hilbident} we have 
$Dh_Z(d+1)>0$. Thus by Proposition \ref{nonincr} and by point iii) of  Lemma \ref{rem:triv} we get 
$Dh_Z(j)>0$ for $j=0,\dots,d+1$. Hence $\sum_jDh_Z(j)\geq d+2$, which contradicts point vii) of  Lemma \ref{rem:triv}.
\end{proof}

An easy extension of Theorem \ref{linbound} is given by the following result. 

\begin{thm0}\label{suplinbound} Assume that a tensor $T\in\Pj(Sym^d(\C^{n+1}))$ has a decomposition $A$ with $\ell(A)\leq (d+n)/2$,
such that $\langle A\rangle=\Pj^n$. Then $T$ has rank $\ell(A)$ and it is identifiable.
\end{thm0}
\begin{proof} Assume on the contrary that $T$ has a second decomposition $B$ with $\ell(B)\leq \ell(A)$, and
take the union $Z=A\cup B$. Then $\ell(Z)\leq 2\ell(A)\leq d+n$. By Lemma \ref{Hilbident} we have 
$Dh_Z(d+1)>0$. Thus by Proposition \ref{nonincr} and by point iii) of  Lemma \ref{rem:triv} we get 
$Dh_Z(j)>0$ for $j=0,\dots,d+1$. By Example \ref{indepcond} and by Proposition \ref{inclu} we get $h_Z(1)=n+1$, so that  $Dh_Z(1)=n$.
Hence $\sum_jDh_Z(j)\geq d+n+1$, which contradicts point vii) of  Lemma \ref{rem:triv}.
\end{proof}

A tensor $T\in\Pj(Sym^d(\C^{n+1}))$ is \emph{concise} if there exist no linear subspaces $W\subset \C^{n+1}$, 
of codimension $1$, such that $T$ belongs to $\Pj(Sym^d(W))$. 

The previous statement implies that when $T$ is concise and it has a decomposition of cardinality $\leq (d+n)/2$, then $T$ is identifiable.
\smallskip

To go further, we may assume some restrictions on the geometry of a decomposition $A$ of $T$. 

\begin{lemma0} \label{linBGM} Let $Z\subset \Pj^n$ be a finite set and assume that for some $j\geq 1$: $Dh_Z(j+1)=Dh_Z(j)=1$.
Then $Z$ contains an aligned  subset $Z'$ of cardinality $\ell(Z')=j+2$, and $Dh_Z(i)=Dh_{Z'}(i)$ for all $i\geq j$.
\end{lemma0}
\begin{proof} See Lemma 2 of \cite{BernGimiIda11}.
\end{proof}

The following result  gives a further extension of Theorem \ref{linbound} (compare with Theorem 2 of  \cite{BallBern12a}).

\begin{prop0}\label{nocoll}
Fix a form $T\in\Pj(Sym^d(\C^{n+1}))$ and a minimal decomposition  $A\subset\Pj^n$  of $T$. Assume that 
 $\ell(A)\leq d$ and $A$ does not contain an aligned subset of cardinality $d/2$. Then $T$ has rank $\ell(A)$ and it is identifiable.
\end{prop0}
\begin{proof} Assume there exists another decomposition $B$ of $T$ with $\ell(B)\leq d$ and call $Z$ the union $Z=A\cup B$.
Then $\ell(Z)\leq 2d$, moreover, by Lemma \ref{Hilbident}, $Dh_Z(d+1)>0$, which implies $Dh_Z(d)>0$. By Example \ref{Exh(1)} we
get that $h_A(1)=2$, hence also $h_Z(1)=2$, by
 Now assume $Dh_Z(d)\geq 2$. Then $Dh_Z(j)\geq 2$ for $j=1,\dots,d$, by Proposition \ref{nonincr}, so that $\sum_jDh_Z(j)\geq 2d+2$, 
which contradicts point vii)  of  Lemma \ref{rem:triv}.  Then for some $j\geq 1$, $j\leq d$, we have $Dh_Z(j)< 2$.
By Proposition \ref{nonincr} again, this implies $Dh_Z(d) =Dh_Z(d+1)=1$. Hence by
lemma \ref{linBGM}, $Z$ contains an aligned subset $Z'$ with $\ell(Z')\geq d+2$, and $Dh_Z(i)=Dh_{Z'}(i)$ for $i>d$. 
Since $Z'$ cannot contain $A$, then there exists a proper subset $A'\subset A$ and a subset $B'\subset B$ such that $Z'=A'\cup B'$.
Shrinking $B'$, if necessary, we may assume that $B'\cap A=\emptyset$, so that also $A'\cap B'=\emptyset$. Then by \eqref{inters}:
\begin{multline*} \dim(\langle \nu_d(A')\rangle\cap\langle\nu_d(B')\rangle) = \ell(Z')-h_{Z'}(d)-1= \sum_{i>d}h_{Z'}(i) =\\
= \sum_{i>d}h_Z(i) = \dim(\langle \nu_d(A)\rangle\cap\langle\nu_d(B_0)\rangle),
\end{multline*}
where $B_0=B\setminus A$. Thus:
$$\langle \nu_d(A')\rangle\cap\langle\nu_d(B')\rangle = \langle \nu_d(A)\rangle\cap\langle\nu_d(B_0)\rangle
$$ 
Since, as in the proof of Lemma \ref{inters},  the intersection  $\langle \nu_d(A)\rangle\cap\langle\nu_d(B)\rangle$ is spanned by
$\nu_d(A\cap B)$ and $\langle \nu_d(A)\rangle\cap\langle\nu_d(B_0)\rangle$, it follows that $T$ belongs to the span
of $\nu_d((A\cap B)\cup A')$. The minimality of $A$ implies $A=(A\cap B)\cup A'$, so the points of $A$ which are
not contained in $B$ are aligned. By assumption $\ell(A')\leq d/2$ and $\ell(A')+\ell(B')=\ell(Z')\geq d+2$, it follows
that $\ell(B')\geq 2+d/2$. Thus $\ell(A\cap B)\leq \ell(B)-2-d/2\leq \ell(A)-2-d/2$. Then 
$$\ell(A)\leq \ell(A')+\ell (A\cap B)\leq d/2+\ell(A)-2-d/2=\ell(A)-2,$$
a contradiction. 
\end{proof}

In order to go further in the study of the identifiability of symmetric tensors, one needs a refinement of lemma
\ref{linBGM}. The refinement is provided by the following, strong result of Bigatti, Geramita and Migliore
(for the case $n=2$ the result has been proved by Davis).

\begin{thm0}\label{BGM} Let $Z\subset\Pj^n$ be a finite set. Assume
that for some $s\leq j$,  $Dh_Z(j)=Dh_Z(j+1)=s$ . Then there exists a reduced curve $C$ of degree $s$
such that, setting $Z'=Z\cap C$ and $Z''=Z\setminus Z'$:
\begin{enumerate}
\item for $i\geq j-1$, $h_{Z'}(i)=h_Z(i)-\ell(Z'')$;
\item for $i\leq j$, $h_Z(i)=h_C(i)$;
\item $Dh_{Z'}(i)=\begin{cases} Dh_C(i) \mbox{ for } i\leq j+1; \\ Dh_Z(i) \mbox{ for } i\geq j.\end{cases}$
\end{enumerate}

In particular, $Dh_{Z'}(i)=s$ for $s\leq i\leq j+1$.

For $n=2$, i.e. when $Z\subset\Pj^2$, we also have:
$$ h_{Z''}(j-1)=\ell(Z'') \quad\mbox{ and }\quad Dh_{Z''}(i)=Dh_Z(i+s)-s  \mbox{ for } i+s\leq j.$$
\end{thm0}
\begin{proof} See Theorem 3.6 of \cite{BigaGerMig94}, and \cite{Davis85} for the case $n=2$.
\end{proof}

Thanks to Theorem \ref{BGM}, for the case $n=2$ one can prove an extension of Proposition \ref{nocoll}:

\begin{thm0} Fix a a form $T\in\Pj(Sym^d(\C^3))$ and a minimal decomposition  $A\subset\Pj^n$  of $T$. Assume that 
 for all $j$ the Kruskal rank of  $v_j(A)$ is maximal, i.e. it is equal to the minimum between  $\ell(A)$ and $\binom{j+2}2$. If
 $$\ell(A)<  \frac{d^2+d}8,$$
 then $T$ has rank $\ell(A)$ and it is identifiable.
\end{thm0}
\begin{proof}
See Theorem 1.4 of \cite{BallC12}, in which the general uniform position (GUP) assumption
is equivalent to the condition that the Kruskal rank of $v_j(A)$ is maximal for all $j$. 
\end{proof}

One aspect of the study of decomposition which has not been developed appropriately derives from the 
observation that Sylvester Theorem \ref{Sy} can be sharpened as follows.

\begin{thm0}\label{Sy2} Assume $n=1$. Assume that $T\in\Pj(Sym^d(\C^{2}))$ has a minimal decomposition $A$ with
$\ell(A)<d+1$. Then for any other minimal decomposition $B$ of $T$ one has $\ell(A)+\ell(B)\geq d+2$. 
\end{thm0}
\begin{proof} 
Assume on the contrary that $T$ has a second decomposition $B$ with $\ell(B)+ \ell(A)\leq d+1$, and
take the union $Z=A\cup B$. Then $\ell(Z)\leq d+1$. Then we conclude as in the
proof of Theorem \ref{linbound}.
\end{proof}

\begin{rem0} {\rm One can prove a statement similar to Theorem \ref{suplinbound} under the  assumption that $\langle A\rangle
=\Pj^n$. Namely in this case for any other minimal decomposition $B$ of $T$ one has $\ell(A)+\ell(B)\geq d+n$. 
Details are left to the reader.}
\end{rem0}

\section{Kruskal's criterion and Terracini's criterion}

The most famous and most used criterion for detecting the identifiability of a given tensor was proved
by Kruskal in 1977 (see \cite{Kruskal77}). Kruskal's criterion was originally proved for 3way,
non necessarily symmetric, tensors. The application to symmetric tensors 
of any size is described e.g. in \cite{COttVan17b}. We recall the result here, rephrased 
in terms of the geometric language. 

\begin{thm0}\label{kruthm} {\bf Reshaped Kruskal's criterion}. 
Let $T\in Sym^d(\C^{n+1})$ and let $A\subset \Pj^n$ be a minimal
decomposition of $T$. Fix a partition $d=a+b+c$, with $0<a\leq b\leq c$. Write
$k_a,k_b,k_c$ for the Kruskal ranks of $v_a(A)$, $v_b(A)$, $v_c(A)$ respectively. If
$$\ell(A)\leq \frac {k_a+k_b+k_c-2}2$$
then $T$ has rank $\ell(A)$ and it is identifiable.
\end{thm0}

Of course the efficiency of the previous criterion depends on the choice of the partition. One
should observe that computing the Kruskal ranks can be  demanding, for large values of $d$, unless
the  coordinates matrices of $v_a(A)$, $v_b(A)$, $v_c(A)$ have full rank. For that reason,
and also for widening the range in which Kruskal's criterion applies, it is usually convenient
to us a maximally unbalanced partition

\begin{exa0} {\rm Consider the case $d=4$. The unique partition  is $a=b=1, c=2$.

If $2\leq\ell(A)\leq n+1$, in the most favorable case in which $k_a=k_b=k_c=\ell(A)$,
then the condition   $\ell(A)\leq(k_a+k_b+k_c-2)/2$ is automatically satisfied and Kruskal's criterion
applies.

If $n+1<\ell(A)\leq \binom{n+2}2$, then the most favorable case is $k_a=k_b=n+1$ and
$k_c=\ell(A)$. In this situation $\ell(A)\leq(k_a+k_b+k_c-2)/2$ is equivalent to $\ell(A)\leq 2n$.

So, one cannot hope to apply directly Kruskal's criterion,for $d=4$, as soon as $\ell(A)>2n$.
}\end{exa0} 

A direct improvement of Kruskal's criterion is impossible, unless one adds some extra test on the tensor $T$.
Namely Kruskal's criterion (even in its reshaped version) is known to be sharp, in its maximal range.

\begin{thm0} For any $n,d,a,b,c$ there exists a tensor $T\in Sym^d(\C^{n+1})$ with a minimal decomposition 
 $A$ such that the Kruskal's ranks $k_a,k_b,k_c$ are maximal (i.e. $k_a=\min\{\ell(A),\binom{n+a}a$,
 and a similar equality holds for $b$ and $c$), with
 $$\ell(A)=\frac{k_a+k_b+k_c}2$$
 and such that $T$ is not identifiable.
\end{thm0}
\begin{proof} The proof is essentially due to Derksen (\cite{Derksen13}), who proved the result
in the non symmetric case. Remark 1.1 of \cite{AngeCVan} contains the observation that, when $T$ is
symmetric, then Derksen's construction provides several \emph{symmetric} decomposions of $T$.
\end{proof}

Thus, given a decomposition $A$ of  a fixed symmetric tensor $T$, one can test the identifiability
(and the rank) of $T$ by computing the Kruskal ranks $k_a$ of the images of $A$ in suitable
Veronese embeddings, hoping to obtain $k_a+k_b+k_c\geq 2\ell(A)+2$. If the inequality holds,
Kruskal's theorem guarantees the identifiability of $T$.

Typically, the reshaped Kruskal's criterion works for small values of $\ell(A)$. To study the
identifiability of tensors in a wider range, one needs to add some new test for $T$.

An example of a test that, together with Kruskal's test, can provide an affirmative 
answer for the identifiability of $T$, is provided by an observation which comes out from the Terracini's
description of the tangent space to the set of tensors of fixed rank.
\medskip

 In the space $\Pj(Sym^d(\C^{n+1}))$, call $\Sigma_r$ the set of tensors of rank $r$. 
 
 For small values of $r$, i.e. for $r(n+1)\leq \binom{n+d}d$, $\Sigma_r$
 is locally closed in the Zariski topology, i.e. it is an open subset of a projective subvariety
 (the $r$-th \emph{secant variety} of the Veronese image $v_d(\Pj^n)$).
 
Consider the symmetric product $(\Pj^n)^{(r)}$. In the product 
$$\Pj(Sym^d(\C^{n+1}))\times (\Pj^n)^{(r)}$$ 
consider the subvariety $A\Sigma_r$ of pairs $(T,[\{P_1,\dots,P_r\}]) $ such that 
the set $A=\{P_1,\dots,P_r\}$ is mapped by $v_d$ to 
a finite set which spans a subspace of dimension $r-1$ in $\Pj(Sym^d(\C^{n+1}))$
(i.e. $v_d(A)$ is linearly independent) and $T$ belongs to the span of $v_d(A)$.

The set $A\Sigma_r$, which is a quasi-projective variety, is called the
\emph{abstract secant variety} of  $v_d(\Pj^n)$. The projection to
the first factor maps $A\Sigma_r$ surjectively to $\Sigma_r$. 

\begin{defn0} {\rm Define the $r$-th \emph{secant map} $s_r$ as the map
projection to the first factor
$$ s_r: A\Sigma_r\to  \Pj(Sym^d(\C^{n+1})).$$}
\end{defn0}

The image of the secant map is $\Sigma_r$. The inverse image of a tensor $T$
of rank $r$ in the secant map is the set of decompositions of $T$.
\smallskip

Since $v_d(\Pj^n)$ is a smooth variety, then $(\Pj^n)^{(r)}$ is smooth,
outside the diagonals. Thus also $A\Sigma_r$, which is a $\Pj^{r-1}$ bundle over 
a subset of $(\Pj^n)^{(r)}$ which does not meet the diagonals, is smooth. 

\begin{defn0} \label{computerr}{\rm The tangent space to 
$A\Sigma_r$ at a point $(T,[\{P_1,\dots,P_r\}]) $ maps, in the differential of $s_r$,
to the space $\mathcal T$ in $\Pj^N=\Pj(Sym^d(\C^{n+1}))$ spanned by the tangent spaces 
to $v_d(\Pj^n)$ at the points $v_d(P_1),\dots,v_d(P_n)$.
We call this space the \emph{Terracini space} of the decomposition 
$A=\{P_1,\dots,P_r\}$ of $T$.}
\end{defn0}

The name of \emph{Terracini space} comes from the celebrated Terracini's Lemma,
which says that, for a general choice of $T\in \Sigma_r$ and for $r\leq N$, the Terracini space is
the tangent space to $\Sigma_r$ at $T$. Thus, a computation of the dimension 
of the Terracini space at a general point corresponds to compute the dimension
of the set $\Sigma_r$ of tensors of rank $r\leq N$.  

\begin{rem0}\label{remdim} {\rm The dimension of the Terracini space $\mathcal T$ is naturally bounded:
$$ \dim(\mathcal T)\leq (n+1)r-1,$$
and the equality means that the tangent spaces to $v_d(\Pj^n)$ at the points $v_d(P_i)$'s
are linearly independent.

Since $A\Sigma_r$ is a $\Pj^{r-1}$ bundle over a quasi-projective variety of dimension $nr$,
then $(n+1)r-1=\dim(A\Sigma_r)$. It follows that the dimension of the Terracini space equals
$(n+1)r-1$ when the differential of $s_r$ has maximal rank.
}\end{rem0}

\begin{rem0} {\rm The decomposition $A$ of  $T\in\Pj(Sym^d(\C^{n+1}))$ corresponds to
the datum of $r$ linear forms $L_1,\dots,L_r$ in the polynomial ring
$R=\C[x_0,\dots,x_n]$.

The Terracini space can be naturally identified with the degree $d$ homogeneous piece
of the ideal in $R$ spanned by 
$$L_1^{d-1}m+ \dots +L_r^{d-1}m,$$
where $m$ is the ideal generated by the variables.

It follows that the computation of the dimension of the Terracini space at a decomposition
of $T$ is a straightforward application of simple algorithm of \emph{linear} algebra.

We refer to the book \cite{IK} for the (elementary) proof of this statement.
}\end{rem0}

The use of the Terracini space in the computation of the identifiability of
a form $T$ is meaningful in the following situation.

\begin{prop0}\label{fam} Let $A$ be a decomposition of $T$ of length $r$ and assume that there exists a 
non trivial family $A_t$ of decompositions of $T$, such that $A_0=A$. Then the Terracini space 
of $A$ has dimension strictly smaller than $(n+1)r-1$.
\end{prop0}
\begin{proof} $A_t$ determines a positive dimensional subvariety $W$ in the fiber of $s_r$ over $T$.
Thus, there exists a tangent vector to $A\Sigma_r$ at $(T,[A])$, where $[A]$ is the point of the 
symmetric product corresponding to $A$, which is killed by the differential of $s_r$ at $(T,[A])$.
Then use remark \ref{remdim}.
\end{proof}

Now, we can introduce our strategy in finding criteria for the identifiability of symmetric tensors,
which works in a range slightly wider than  the Kruskal's one.

If we can prove that tensors $T$ which are non identifiable must have  a 
positive dimensional family of different decompositions, containing the given decomposition $A$,
then  we can check the identifiability of $T$ by computing the dimension of the Terracini space.

The fact that non identifiable tensors have indeed a positive dimensional family of different decompositions,
is false in general. It turns out, however, that this fact holds in some cases, especially when we are outside the
Kruskal's numerical range, but very close to it.
\smallskip

A way to produce positive dimensional family of different decompositions is explained in the following:

\begin{prop0}\label{curvprop} Assume that a decomposition $A$ of length $r$ of $T$ is contained in a projective curve $C\subset\Pj^n$
which is mapped by $v_d$ to a space $\Pj^m$, with $m<2r-1$. Then there exists
positive dimensional family of different decompositions $A_t$ of $T$, such that $A_0=A$.
\end{prop0}
\begin{proof} $T$ belongs to the span of $v_d(A)$, which is contained in the span of $v_d(C)$,
which is contained in $\Pj^m$. The condition $m<2r+1$ implies that there is a positive dimensional
family of subsets $A_t\subset C$ such that $T\in\langle v_d(A_t)\rangle$. Namely, the abstract $r$ secant
variety $A\Sigma^C_r$ of $C$ has dimension $2r-1$, thus all the components of the fibers of the 
map $A\Sigma^C_r\to\Pj^m$ are positive dimensional.
\end{proof}

Now we can mix together the analysis of the Hilbert function, the Cayley-Bacharach condition and
the computation of the dimension of the Terracini space, to produce a criterion for the identifiability of $T$.

\begin{thm0}\label{quart} {\bf (See \cite{AngeCVan})}. Let $T$ be a quartic form in $n+1$ variables,
and consider a decomposition $A$ of $T$ of length $2n+1$.

Assume that:
\begin{itemize}
\item[a)] the Kruskal rank of $A$ is $n+1$;
\item[b)] the Terracini space at $A$ has (the maximal) dimension $(2n+1)(n+1)-1$.
\end{itemize}

Then $T$ has rank $2n+1$ and it is identifiable.
 \end{thm0}

Notice that conditions a) and b) are expected to hold for a general quartic, i.e. outside a proper 
Zariski closed subset  (of measure $0$) in the space of quartics. Thus the previous theorem provides a criterion 
to prove the identifiability of $T$, except for very special tensors.

\begin{proof} We give a sketch of the proof.

First notice that, by Proposition \ref{v2}, the set $v_2(A)$ is linearly independent, i.e. it has Kruskal rank $2n+1$-

Call $B$ a different decomposition of length $\leq 2n+1$ for $T$, which we want to exclude. Call $Z$ the union
$Z=A\cup B$ and consider the Hilbert function of $Z$. 

First step is to prove that  $Z$ has the Cayley-Bacharach property $CB(4)$. This is almost clear when
$A\cap B=\emptyset$, while if $A\cap B\neq \emptyset$ the claim follows from Kruskal's theorem.

Next, since $Z$ has the property $CB(4)$, by Theorem \ref{GKRext}
it follows soon that $Dh_Z(3)+Dh_Z(4)+Dh_Z(5)\geq h_A(2)=2n+1$,
so that $h_Z(2)=h_A(2)=2n+1$. Then one invokes the following extension of the
classical Castelnuovo's Lemma:

\begin{lemma0}\label{Cnew} {\bf(See \cite{AngeCVan}, Lemma 5.4) }
 Let $Z$ be a set of $r\geq 2n+3$ points in $\Pj^n$ which impose 
 at most $2n+1$   conditions to quadrics. Assume that $Z$ has a 
 subset $Z'$ of $2n+1$  points in LGP. Then the entire $Z$ is 
 in LGP and it is contained in an irreducible  rational normal curve.
\end{lemma0}

The classical formulation of Castelnuovo's lemma required that the whole set $Z$ is in LGP,
which we cannot assume in our setting, because we only know the position of $A$, which contains 
$2n+1$ points of $Z$, while we have no control of the points of $B$. Fortunately, the extension
matches exactly our requirements. Now we can turn back to the proof of the Theorem.

Since $Z$ has a subset, namely $A$, which is in LGP, then it follows that $Z$, hence also $A$,
sits in a rational normal curve $C$ of $\Pj^n$. The image of $C$ in the Veronese map
$v_4$ spans a $\Pj^{4n}$. Hence the claim follows by Proposition \ref{curvprop}. 
\end{proof}

\section{A new result on the decomposition of tensors}

In this section we improve slightly Theorem \ref{quart}, by removing the assumption that
the Kruskal rank of $A$ is $n+1$, and replacing it by a numerical assumption on $\ell(A)$.
At a certain point of the proof we will need the cohomological properties of the residue
of a finite set with respect to a hyperplane. This is the unique passage in which some
sophisticated algebraic machinery enters into the proof.

Let $Z$ be a finite set in $\Pj^n$, and let $H$ be a hyperplane. Call $Z_1$ the intersection 
$Z_1=Z\cap H$ and call $Z_2$ the set:
$$ Z_2= Z\setminus Z_1 = Z\setminus (Z\cap H).$$
For obvious reasons, $Z_2$ is called the {\it residue} of $Z$ with respect to $H$. 

If $I_Z,I_{Z_2}$ denote the homogeneous ideals of $Z,Z_2$ respectively, the multiplication
by an equation of $H$ determines an exact sequence of graded modules:
\begin{equation}\label{eqres}  0\to I_{Z_2}(1) \to  I_Z(2) \stackrel{\rho}\longrightarrow I_{Z_1,H}(2) \end{equation}
in which the rightmost ideal $I_{Z_1,H}$ is the homogeneous ideal of $Z_1$ in $H$.

The following result is a straightforward application of the cohomology of maps of sheaves:

\begin{lemma0}\label{reslemma} Assume that $Z_2$ is linearly independent. Then the rightmost
map $\rho$ in sequence \eqref{eqres} is surjective.
\end{lemma0}
\begin{proof} The cokernel of $\rho$ is contained in the cohomology group $H^1(\mathcal I_{Z_2}(1))$,
where $\mathcal I_{Z_2}$ is the ideal sheaf of $Z_2$. Moreover $H^1(\mathcal I_{Z_2}(1))$ vanishes if
$Z_2$ is linearly independent, because in this case the evaluation map $ev(1)$ on $Z_2$ determines a surjective map
$\C^{n+1}\to \C^{\ell(Z_2)}$.
\end{proof}

\begin{rem0}{\rm With the same trick, one can prove the following general statement:

Assume that the residue $Z_2$ of a finite set $Z$, with respect to a hyperplane $H$, is separated by forms
of degree $d-1$. Then any form of degree $d$ in $H$ which contains $Z_1=Z\setminus Z_2$ can be 
lifted to a form of degree $d$ in $\Pj^n$ which contains $Z$.
}
\end{rem0}

As a consequence, in the hypothesis of Lemma \ref{reslemma}, it turns out that every quadric of the
hyperplane $H$ that contains $Z_2$ can be lifted to a quadric of $\Pj^n$ that contains $Z$.

We will need the following, well known remark for linearly independent sets $W$
in a projective space $\Pj^n$:

\begin{lemma0}\label{indq}  Let $W$ be a linearly independent finite set in $\Pj^n$.
Then for any $Q\notin W$, there exists a quadric of $\Pj^n$ containing $W$
and missing $Q$. In other words, the ideal of $W$ is generated by quadrics.
\end{lemma0}
\begin{proof}
The proof is an easy argument of linear algebra. After shrinking $n$ we may
always assume $W=\{P_1,\dots,P_{n+1}\}$. If $Q$ does not belong to the span of any proper subset of $W$, just by taking
two hyperplanes containing two proper subsets we get the claim. Thus, reorder the points
of $W$ so that $P_1,\dots, P_s$ ($s\geq 2$, $s\leq n$) is a minimal subset whose span $L$ contains $Q$. 
Since the points are linearly independent, the span $M$ of $P_1,\dots,P_{s-1},P_{s+1}$ intersects
$L$ in the span of $P_1,\dots,P_{s-1}$, hence by minimality it does not contain $Q$.
Similarly, the span $M'$ of $P_s,P_{s+2},\dots,P_{n+1}$ intersects $L$ only in $P_n$.
The union of a general hyperplane containing $M$ and a general hyperplane containing $M'$
provides a quadric containing $W$ and missing $Q$.
 \end{proof}

Now we are ready to state and proof our result.

\begin{thm0}\label{quartplus} . Let $T$ be a quartic form in $n+1$ variables,
and consider a decomposition $A$ of $T$. Call $k$ the Kruskal rank of $A$ and assume that $\ell(A)\leq 2k-1$.
Assume that the Terracini space at $A$ has (maximal) dimension $(2k-1)(n+1)-1$.

Then $T$ has rank $2k-1$ and it is identifiable.
 \end{thm0}
 
Notice that since $k\leq n+1$, it follows $\ell(A)\leq 2n+1$. Moreover, by Proposition \ref{v2},
we know that $A$ is separated by quadrics, i.e. $v_2(A)$ is linearly independent. This implies immediately that
also $v_4(A)$ is linearly independent.

 Notice also that if $\ell(A)<2k-1$, then $A$ satisfies the hypothesis of the reshaped Kruskal's criterion,
 because in this case 
 $$ \ell(A)\leq \frac{k+k+\ell(A)-2}2,$$
 so that the identifiability of $A$ follows immediately. 
 
 Thus the Theorem produces a new criterion only for $\ell(A)=2k-1$. Hence we assume, in the proof,
 that $\ell(A)=2k-1$.
 
\begin{proof} As in the proof of Theorem \ref{quart}, we will prove that, under the assumptions,
if another decomposition $B$ of cardinality $\ell(B)\leq 2k-1$ exists, then there exists a curve $C$ containing
$A$ and such that $v_4(C)$ spans a space of dimension $\leq 4k-4$, which contradicts the assumption  2).

Of course, we may assume that $A$ spans $\Pj^n$, otherwise we simply decrease $n$. It follows that $2k-1>n$ and
the difference of the Hilbert function of $A$ is:
$$ Dh_A(0)=1,\quad Dh_A(1)=n, \quad Dh_A(2)=2k-2-n.  $$
Assume that a second decomposition $B$ exists. The first step is to prove that $Z=A\cup B$ satisfies the
Cayley Bacharach property $CB(4)$, which holds by following verbatim
the proof of the similar statement in Theorem 6.2 of \cite{AngeCVan}.

It follows then, by  Theorem \ref{GKRext}, that the difference of the Hilbert function
of $Z$ satisfies $Dh_Z(3)+Dh_Z(4)+Dh_Z(5)=2k-1$, so that in particular $\ell(B)=2k-1$,
$A,B$ are disjoint and the difference of the Hilbert function of $Z$ satisfies:
$$ Dh_Z(0)=1,\quad Dh_Z(1)=n, \quad Dh_Z(2)=2k-2-n.  $$
Thus, summing up, one gets $h_Z(2)=h_A(2)$, i.e. all the quadrics that contain $A$ must contain $Z$.

The assumption that $k$ is the Kruskal rank of $A$ means that any subset of $k$ points in $A$
is linearly independent, while there exists a subset of $k+1$ points which generates a
subspace  $\Lambda=\Pj^{k-1}$. After rearranging the points, we may assume that
$P_1,\dots,P_{k+1}$ generate $\Lambda$, $P_{k+1},\dots,P_{k+q}$ are also contained in
$\Lambda$, and $P_{k+q+1},\dots,P_{2k-1}$ are outside $\Lambda$. Notice that we may always assume
that $A$ is non degenerate, thus $k+q<2k-1$. Call $\Lambda'$ the space generated by $P_{k+q+1},\dots,P_{2k-1}$.
Any pair of hyperplanes $H,H'$ which contain $\Lambda,\Lambda'$ respectively, determine a quadric which
contains $A$. It follows that all the points of $B$ are contained either in $\Lambda$ or in $\Lambda'$.

Let $Q$ be a point of $B$ which lies in $\Lambda$. For any subset $W$ of $k-1$ points among $P_1,\dots,P_{k+q}$
consider the hyperplane $L_W$ of $\Lambda$ spanned by $W$. If $Q$ belongs to no hyperplanes $L_W$, then there are
quadrics in $\Lambda$ which contain $P_1,\dots,P_{k+q}$. Thus if $H$ is a general hyperplane containing $\Lambda$ then
there are quadrics in $H$ which contain $P_1,\dots,P_{k+q}$ and miss $Q$. Since the set 
$P_{k+q+1},\dots,P_{2k-1}$ is linearly independent, by our assumption on the Kruskal rank of $A$, then
by Lemma \ref{reslemma} one finds a quadric of $\Pj^n$ which contains $A$ and misses $Q$, contradicting
$h_Z(2)=h_A(2)$. 

Hence, there exists a set $W$ of $k-1$ points among $P_1,\dots,P_{k+q}$ which spans a hyperplane $L_W$
of $\Lambda$ containing $Q$. Since $W$ is linearly independent, by Lemma \ref{indq} one can find a quadric $K$ in 
$L_W$ that contains $W$ and misses $Q$. Since, by our assumption on the Kruskal rank of $A$, also
$\{P_1,\dots,P_{k+q}\}\setminus W$, which contains at most $k$ points, is linearly independent, then
by Lemma \ref{reslemma} we can lift $K$ to a quadric $K'$ of $\Lambda$ which misses $Q$ and contains
$P_1,\dots,P_{k+q}$. As above, $K'$ lifts to a quadric $K''$ which contains $A$ and misses $Q$.
Thus we have a contradiction with $h_Z(2)=h_A(2)$.

It follows that all the points of $B$ belong to $\Lambda'$. In particular, the form $T$ does not involve
all the variables. After choosing carefully the coordinates $x_0,\dots,x_n$ in $\Pj^n$, we may assume that $T$ does not
involve $x_n$. But then, by replacing $x_n$ with $tx_n$ in the points of $A$
(actually in the points of $A\cap\Lambda$), as $t$ varies we get a family of decompositions of $T$
which coincides with $A$ for $t=1$. By Proposition \ref{fam}, this contradicts the assumption that the Terracini 
space has maximal dimension.
 \end{proof}

\begin{rem0} {\rm As in section 6 of \cite{AngeCVan}, one can create an algorithm
that uses Theorem \ref{quartplus} to detect the identifiability of quartics of low rank.
Given a  symmetric  decomposition of length $r$ of a quartic 
$$  T = \sum_{i=1}^r \nu_{4}( P_i ), $$
in the form of the collection of points $A = \{ P_i = [\vect{m}_i] \}_{i=1}^r \subset \Pj^n$, we can apply the following 
algorithm for verifying that the given decomposition of $T$ is identifiable:
\begin{enumerate}
  \item[1)] \emph{Kruskal's test}: compute the Kruskal rank $k$ of  $A$;
 \item[S1.] If $r > 2k-1$, the criterion cannot be applied.
 \item[S2.] If $r < 2k-1$, use the reshaped Kruskal criterion from \cite{COttVan17a}, section 6.2.
 \item[S3.] If $r = 2k-1$, perform the: 
  \item[2)] \emph{Terracini's test}: check that $\dim \langle \Tang{\vect{m}_1}{\nu_4(\C^{n+1})}, \ldots, \Tang{\vect{m}_r}
  {\nu_4(\C^{n+1})} \rangle = (2k-1)(n+1)-1$.
  \end{enumerate}
  
 If all these tests are successful, then $T$ is of rank $r$ and is identifiable.
}\end{rem0}

Notice that the  computation of the Kruskal rank of $A$  turns out to be the heaviest step of the algorithm. 

\section{Final remarks and open problems}

1. We believe that the range in which the non-identifiability of tensors
implies the existence of a positive dimensional family of decompositions (which can
be detected by the computation of the Terracini space) goes beyond the numerical
bounds given in theorems \ref{quart} and \ref{quartplus}.

In order to extend the previous results, however, one needs extensions of
 the basic Castelnuovo's Lemma
\ref{Cnew}. What we would need is to replace the existence of a rational normal curve,
predicted by Lemma \ref{Cnew} for sets of points with special Hilbert functions,
with the existence of other types of curves (elliptic, or even of higher genera), when
the number of points increases.

Similar results are known in some cases (see e.g. \cite{Petrakiev08}, \cite{Ghezzi10}), but not in a form
that can be immediately applied to our situation.

We would like to stimulate further  researches on the geometry of sets of points with
special Hilbert functions, with the final target of an application
to tensor analysis.

\bigskip

2. The geometric methods known so far for the study of the identifiability of  specific tensors,
as the Kruskal's criterion and the extension given in the previous sections, are based on the study
of the geometry of a given decomposition. The idea has a basic bug: once the identifiability
follows from geometric properties of a given decomposition $A$, then it must hold for all
the tensors which lie in the span of $v_d(A)$ (at least those for which $A$ is minimal), regardless of the
 coefficients that are used to produce the  form $T$. Of course, we can expect that 
a similar uniform behavior holds only for small values of the rank $r$. When $r$ increases,
then it is natural to expect that the space $\langle v_d(A)\rangle$ contains both
identifiable and non-identifiable points.

As a consequence, we need criteria for identifiability which are able to distinguish between
different points of the span $\langle v_d(A)\rangle$ of a given decomposition $A$.

We believe that a geometric analysis of  $A$ and of its linked sets of points
can produce geometric criteria which reach much further than the range 
of application of Kruskal's criterion. 

\bigskip

3.  A different approach to the study of the identifiability of tensors is contained in the paper 
\cite{MassaMellaStagliano}. The authors prove that when the space spanned by partial derivatives of the form $T$
(the \emph{catalecticant space}, in the terminology of \cite{IK}) meets the corresponding variety in a finite
set of the expected length $r$, then $r$ is the rank of $T$ and the tensor is identifiable.

The method of partial derivatives has the advantage that it does not need to start with a given decomposition.
On the other hand, for special tensors, it does not describe the geometric situation which yields
the non-uniqueness of the  decomposition. Furthermore, the method relies on the computation of an intersection of 
algebraic varieties, i.e. on methods of computer algebra, which usually cost a lot in terms
of computational complexity.

We believe that a mix of the two methods, which will be the target of a forthcoming 
paper, will produce new, interesting developments in the theory.

\bigskip

4. We wonder if the analysis of tensor decomposition by means of geometric methods,
related with the study of finite sets in projective spaces, can be extended beyond the case of symmetric
tensors. For general tensors, the natural substitute for the Hilbert function is the multiugraded
Hilbert function. Indeed, for general tensors, one has only to consider the first piece of
the multigraded Hilbert function, i.e. the piece bounded by the origin and the multidegree $(1,\dots,1)$.
For this piece of the Hilbert function, which is basically the \emph{Segre function}, in the terminology of 
\cite{CSacchi} and \cite{BallBernCGuardo}, very few is known. For instance, we do not know an analogue
of Lemma \ref{rem:triv}, which lists the most elementary properties.

A study of the Segre function, aimed to an application to tensor analysis, will probably yield several new,
valuable results on the theory.


\begin{thebibliography}{alpha}
\bibitem {AngeCVan} E. Angelini, L. Chiantini, and N. Vannieuwenhoven. \emph{Identifiability beyond Kruskal's
bound for symmetric tensors of degree 4.} Rend. Lincei Matem. Applic. 29 (2018), 465-485.
\bibitem {AlexHir95} J. Alexander and A. Hirschowitz. \emph{Polynomial interpolation in several variables}. J. 
Algebraic Geom. 4 (1995), 201-222.
\bibitem {Ball05a}  E. Ballico. \emph{On the weak non-defectivity of Veronese embeddings of projective spaces}.
Central Eur. J. Math. 3 (2005), 183-187.
\bibitem {BallBern12a}  E. Ballico and A. Bernardi. \emph{Decomposition of homogeneous polynomials with low rank}.
Math. Zeit. 271 (2012), 1141-1149.
\bibitem {BallBernCGuardo} E. Ballico, A. Bernardi, L. Chiantini and E. Guardo. \emph{Bounds on the tensor rank}.
Ann. Mat. Pura Appl. (to appear). ArXiv:1705.02299.
\bibitem {BallC12}  E. Ballico and L. Chiantini. \emph{A criterion for detecting the identifiability of symmetric
tensors of size three}. Diff. Geom. Applic. 30 (2012), 233-237.
\bibitem {BernGimiIda11} A. Bernardi, A. Gimigliano, and M. Id\'a. \emph{Computing symmetric rank for symmetric
tensors}. J. Symbolic Comput. 46 (2011), 34-53.
\bibitem {BucGinenskyLand13}  J. Buczy\'nski, A. Ginensky, and J.M. Landsberg. \emph{Determinantal equations for secant
varieties and the Eisenbud-Koh-Stillman conjecture}. J. London Math. Soc. 88 (2013), 1-24.
\bibitem {BigaGerMig94}  A.M. Bigatti, A.V. Geramita, and J. Migliore. \emph{Geometric consequences of extremal
behavior in a theorem of Macaulay}. Trans. Amer. Math. Soc. 346 (1994), 203-235.
\bibitem {CarlCatalC15}  E. Carlini, M.V. Catalisano, and L. Chiantini. \emph{Progress on the symmetric Strassen
conjecture}. J. Pure Appl. Algebra 219 (2015), 3149-3157.
\bibitem {CCi06}  L. Chiantini and C. Ciliberto. \emph{On the concept of k-secant order of a variety}. J. London
Math. Soc. 73 (2006), 436-454.
\bibitem {COttVan17a}  L. Chiantini, G. Ottaviani, and N. Vannieuwenhoven. \emph{Effective criteria for specific
identifiability of tensors and forms}. SIAM J. Matrix Anal. Appl. 38 (2017), 656-681.
\bibitem {COttVan17b}  L. Chiantini, G. Ottaviani, and N. Vannieuwenhoven. \emph{On generic identifiability of
symmetric tensors of subgeneric rank}. Trans. Amer. Math. Soc. 369 (2017), 4021-4042.
\bibitem {CSacchi} L. Chiantini and D. Sacchi. \emph{Segre functions in multiprojective spaces and tensor analysis}.
In: From Classical to Modern Algebraic Geometry, G. Casnati et. Al. Editors. Trends Hist. Sci. 8,  Birkhauser
(2016), 361-374. 
\bibitem {Davis85}  E. Davis. \emph{Complete intersections of codimension $2$ in $\mathbb P^r$: the Bezout-Jacobi-Segre theorem revisited}. Rend. Seminario Mat. Univ. Politecnico Torino 43 (1985), 333-353.
\bibitem {Derksen13}  H. Dersken. \emph{Kruskal's uniqueness inequality is sharp}. 
Linear Alg. Applic. 438 (2013), 708-712.
\bibitem {Ghezzi10}  L. Ghezzi. \emph{A generalization of the strong Castelnuovo lemma}. 
J. of Algebra 323 (2010), 1018-1035.
\bibitem {GerKreuzerRobbiano93}  A.V. Geramita, M. Kreuzer, and L. Robbiano. \emph{Cayley-Bacharach 
schemes and their canonical modules}. Trans. Amer. Math. Soc. 339 (1993), 443-452.
\bibitem {GaluppiMella}  F. Galuppi and M. Mella. \emph{Identifiability of homogeneous polynomials and 
Cremona transformations}. preprint arXiv:1606.06895.
\bibitem {IK} A. Iarrobino and V. Kanev. {\bf Power Sums, Gorenstein Algebras, and Determinantal
Loci}. Volume 1721 of Lecture Notes in Mathematics. Springer, Berlin, New York NY, 1999.
\bibitem {Kruskal77}  J.B. Kruskal. \emph{Three-way arrays: rank and uniqueness of trilinear decompositions, with
application to arithmetic complexity and statistics}. Linear Algebra Appl. 18 (1977), 95-138.
\bibitem {MassaMellaStagliano}  A. Massarenti, M. Mella and G. Staglian\'o. \emph{Effective identifiability criteria for tensors and polynomials}. J. Symbolic Comput. 87 (2018), 227-237.
\bibitem {Mella06}  M. Mella. \emph{Singularities of linear systems and the Waring problem}. Trans. Amer. Math.
Soc. 358 (2006), 5523-5538.
\bibitem {Migliore}  J. Migliore. {\bf Introduction to Liaison Theory and Deficiency Modules}. Volume 165 of
Progress in Mathematics. Birkhauser, Basel, 1998.
\bibitem {Petrakiev08} I. Petrakiev. \emph{A step in Castelnuovo theory via Groebner bases}. J. Reine Angew. Math.,
619 (2008), 49-73.
\bibitem {RaoLiZhang18}  W. Rao, D. Li, and J.Q. Zhang. \emph{A tensor-based approach to L-shaped arrays processing with enhanced degrees of freedom}. IEEE Signal Proc. Lett. 25 (2018), 1-5.
\bibitem {Shitovb} Y. Shitov. \emph{A counterexample to {S}trassen's direct sum conjecture}. Preprint 
arXiv:1712.08660.
\bibitem {Sylvester86}   J.J. Sylvester. \emph{Sur une extension d`un th\'eor\'eme de Clebsch relatif aux courbes du
quatri\'e me degr\'e}. C. R. Math. Acad. Sci. Paris 102 (1886), 1532-1534.
\end{thebibliography}
\end{document}